\RequirePackage[reqno]{amsmath}
\documentclass[a4paper,12pt]{amsart}

\usepackage{empheq}
\numberwithin{equation}{section}
\usepackage{mathrsfs}
\usepackage[active]{srcltx}
\usepackage{longtable}
\usepackage{empheq}

\usepackage{color}
\usepackage[unicode]{hyperref}
\hypersetup{
	colorlinks = true,%
	citecolor = [rgb]{0.0,0.0,0.9},
	filecolor=black,%
	linkcolor = [rgb]{0.65,0.0,0.0},%
	anchorcolor = red,
	pagecolor = red,
	urlcolor= [rgb]{0.65,0.0,0.0}
}

\newtheorem{example}{ Example}[section]
\newtheorem{proposition}{Proposition}[section]
\newtheorem{theorem}{Theorem}[section]

\newtheorem{remark}{Remark}[section]
\numberwithin{equation}{section}

\begin{document}

\title[Poincar\'e's lemma on some non-Euclidean structures]{Poincar\'e's lemma on some non-Euclidean structures}

	\author{Alexandru Krist\'aly}
\address{Department of Economics, Babe\c s-Bolyai University,   400591 Cluj-Napoca,
	Romania \& Institute of Applied Mathematics, \'Obuda
	University,
	1034 Budapest, Hungary}
\email{alexandrukristaly@yahoo.com; kristaly.alexandru@nik.uni-obuda.hu}

\begin{abstract}
In this paper we prove the Poincar\'e  lemma on some $n$-dimensional corank 1 sub-Riemannian structures, formulating the $\frac{(n-1)n(n^2+3n-2)}{8}$ necessarily and sufficiently 'curl-vanishing' compatibility conditions. In particular, this result solves partially an open problem formulated by Calin and Chang.  Our proof is based on a Poincar\'e lemma stated on  Riemannian manifolds and a suitable Ces\`aro-Volterra path integral formula established in local coordinates. As a byproduct, a Saint-Venant lemma is also provided on generic Riemannian manifolds. 
Some examples are presented on the hyperbolic space and Carnot/Heisenberg groups. 
\end{abstract}

	\dedicatory{Dedicated to 
	Professor Philippe G. Ciarlet on the occasion of his 80th birthday}

\maketitle

\vspace{0.2cm}
\section{Introduction and Main result}
Let $\Omega\subseteq \mathbb R^n$ be an open, simply connected set, and $\textit{\textbf{a}}=(a_i)\in \textit{\textbf{C}}^1(\Omega;\mathbb R^n)$, $n\geq 2$. The classical Poincar\'e lemma says that there exists $u\in C^2(\Omega)$ with $${{\bf\nabla}}u=\textit{\textbf{a}}\ {\rm in}\ \Omega,$$ if and only if $\operatorname{\textbf{curl}} \textit{\textbf{a}}=\textbf{0}$ in $\textit{\textbf{C}}(\Omega;\mathbb R^n)$, i.e., 
$${\partial_{x_i} a_j}={\partial_{x_j} a_i}\ {\rm in}\ C(\Omega)\ \ {\rm for\ every}\ i,j=1,...,n.$$ 
Here, as usual, ${{\bf\nabla}}u=({\partial_{x_i} u})\in \textit{\textbf{C}}^1(\Omega;\mathbb R^n)$.  For a weak version of the  Poincar\'e lemma (e.g. in $L^2(\Omega)$) and its equivalent formulation in terms of fundamental results in the theory of PDEs, we refer the reader to Amrouche,  Ciarlet and  Mardare \cite{Amrouche-Ciarlet-Mardare, Amrouche-Ciarlet-Mardare-2} and to the comprehensive monograph by Ciarlet \cite[Chapter 6]{Ciarlet-2013}.  

Very recently, Poincar\'e's lemma has been extended to some specific \textit{low}-dimensional sub-Riemannian structures with rank 2 distributions; e.g., the first Heisenberg group $\mathbb H^1$, Engel-type manifolds, Grushin and Martinet type distributions, and the sub-Riemannian 3-dimensional sphere $\mathbb S^3$, see Calin,  Chang and  Eastwood \cite{Calin-E-1, Calin-E-2} and Calin, Chang and Hu \cite{Calin-1}-\cite{Calin-3}. In the sub-Riemannian setting, the number of equations in the system which is going to be solved is strictly less than the space dimension. Accordingly, the solvability of such gradient-type systems deeply depend on the Lie  bracket generating properties of the sub-Riemannian distributions, and it turns out that the 'curl-vanishing' characterization of the solvability of the sub-Riemannian  system becomes a system of PDEs containing higher-order derivatives.  In order to visualize this phenomenon, we consider the first Heisenberg group  $\mathbb H^1=\mathbb C\times \mathbb R$ endowed with its usual group operation and left-invariant vector fields $X_1=\partial_{x_1}-2x_2\partial_{x_3}$ and $X_2=\partial_{x_2}+2x_1\partial_{x_3}$. The sub-Riemannian system 
\begin{equation}\label{sub-1}
X_1u=a_1,\ X_2u=a_2
\end{equation}
is solvable in $\mathcal F(\mathbb H^1)$(=the space of smooth functions on $\mathbb H^1$) for $\textit{\textbf{a}}=(a_1,a_2)\in \textit{\textbf{C}}^1(\mathbb H^1;\mathbb R^2)$ if and only if
\begin{equation}\label{sub-2}
X_1^2 a_2=(X_1X_2+[X_1,X_2])a_1,\ X_2^2 a_1=(X_2X_1+[X_2,X_1])a_2,
\end{equation}
see e.g. Calin and Chang \cite[Theorem 2.9.8]{Calin-book}.
In addition, the solution $u$ of (\ref{sub-1}) can be given the work done by the force vector field $X=a_1X_1+a_2X_2$ along any horizontal curve starting from  ${\bf 0}\in \mathbb H^1$, called also as the Ces\`aro-Volterra horizontal path integral.

The purpose of our paper is to prove Poincar\'e lemmas on some sub-Riemannian structures of arbitrary dimension with corank 1 distribution, including for instance step-two Carnot groups with not necessarily trivial kernel. In the sequel, we present our main result (see Section \ref{section-3} for the notions used below).  

Let $(M,\mathcal D,g)$ be an $(n+1)$-dimensional sub-Riemannian manifold $(n\geq 2)$, and consider the distribution $\mathcal D$ in a given local coordinate system $(x_i)_{i=1,...,n+1}$ containing vector fields of the form 
\begin{equation}\label{vector-fields}
X_i=\partial_{x_i}+A_i\partial_{ x_{n+1}}, \ i=1,...,n,
\end{equation}
where $A_i:M\to \mathbb R$ are smooth functions depending only on the first $n$ variables, i.e., $A_i=A_i(x_1,...,x_n)$. We assume that
\begin{equation}\label{A-feltetel}
\partial_{x_i}A_j-\partial_{x_j}A_i= c_{ij}\in \mathbb R\ \ {\rm for\ every}\  i,j=1,...,n,
\end{equation}
and $$I_0=\{(i,j): c_{ij}\neq 0\}\neq \emptyset.$$ Due to the latter assumptions, the rank $n$  distribution $\mathcal D$ is nonholonomic on $M$, since 
\begin{equation}\label{Xi-k-c-ik}
[X_i,X_j]=c_{ij}\partial_{x_{n+1}}\ \ {\rm for\ every}\  i,j=1,...,n.
\end{equation}
Given $\textit{\textbf{a}}\in \Gamma(\mathcal D)$(=the set of horizontal vector fields on $M$), we are going to study the solvability of the system 
\begin{equation}\label{sub-riem-equation}
{\bf\nabla}_H u=\textit{\textbf{a}}\ \ \ {\rm in}\ M,
\end{equation} 
where $u\in \mathcal F(M)$ and ${\bf\nabla}_H$ denotes the horizontal gradient. Our main result, the Poincar\'e lemma on sub-Riemannian manifolds, reads as follows: 

\begin{theorem}\label{fotetel}
Let	$(M,\mathcal D,g)$ be an $(n+1)$-dimensional simply connected sub-Riemannian manifold $(n\geq 2)$, where the distribution $\mathcal D$ is given by the vector fields in {\rm (\ref{vector-fields})} with functions $A_i$ depending only on the first $n$ variables, verifying {\rm (\ref{A-feltetel})} and $I_0\neq \emptyset.$

 Given $\textit{\textbf{a}}\in \Gamma(\mathcal D)$, the sub-Riemannian system {\rm (\ref{sub-riem-equation})} has a solution $u\in \mathcal F(M)$ if and only if
 		\begin{empheq}[left=\empheqlbrace]{align}
 		c_{kl}(X_i\tilde a_j-X_j\tilde a_i)&=c_{ij}(X_k\tilde a_l-X_l\tilde a_k)\ \ {\rm for\ every}\ \ i,j,k,l=1,...,n; \label{full-jellemzes}\\
 X_kX_i\tilde a_j-X_kX_j\tilde a_i&=[X_i,X_j]\tilde a_k \ \ \ \ \ \ \ \ \  \ {\rm for\ every}\ \ i,j,k=1,...,n,\label{full-jellemzes-2}
 	\end{empheq}
 		where $\textit{\textbf{a}}= a_iX_i$ and $\tilde a_j= g_{ij}a_i$ $($the summations being from $1$ to $n$$),$  and $(g_{ij})$ are the components of $g$ with respect to the distribution $\mathcal D$.
 Moreover, if $x_0\in M,$ the solution $u:M\to \mathbb R$ for the system {\rm (\ref{sub-riem-equation})} can be obtained by 
 \begin{equation}\label{C-V-1}
 u(x)=c_0+\int_0^1 g(\textit{\textbf{a}}(\gamma(t)),\dot \gamma(t))dt,\ x\in M, 
 \end{equation}
where $c_0=u(x_0)\in \mathbb R$ and $\gamma:[0,1]\to M$ is any horizontal curve joining $x_0$ with $x.$ 
\end{theorem}

Some remarks are in order. 

\begin{remark}\rm 
	(a) Although (\ref{full-jellemzes}) and (\ref{full-jellemzes-2}) contain $n^4$ and $n^3$ conditions, a simple combinatorial reasoning shows that it is enough to verify at most $s_n=\frac{(n-2)(n-1)n(n+1)}{8}$ and $s_n'=\frac{(n-1)n^2}{2}$ conditions, respectively. 
	Thus, the number of compatibility conditions is $s_n+s_n'=\frac{(n-1)n(n^2+3n-2)}{8}.$
	
	(b) Theorem \ref{fotetel} provides an answer to the open question of Calin and Chang \cite[p. 55]{Calin-book} whenever the sub-Riemannian manifold with arbitrarily dimension has corank 1 distribution. We note that the existing results in the literature solve the system (\ref{sub-riem-equation}) only for two components, i.e., the distributions contain two vector fields. In particular, if $M=\mathbb H^1$ is the first Heisenberg group, the solvability of the system (\ref{sub-1}) can be recovered by Theorem \ref{fotetel}; indeed, in this particular case, $n=2$, $\mathcal D=\{X_1,X_2\}$ and $g_{ij}=\delta_{ij}$. Moreover, $A_1=-2x_2$, $A_2=2x_1$; thus $c_{12}=-c_{21}=4$ and $c_{11}=c_{22}=0$ in  (\ref{A-feltetel}). Notice that the first-ordered relations in (\ref{full-jellemzes}) are trivially satisfied (supported also by the fact that $s_2=0$, thus nothing should be checked), while the second-ordered ones (\ref{full-jellemzes-2}) reduce precisely to (\ref{sub-2}), containing  $s_2'=2$ conditions. In higher-dimensional Heisenberg groups $\mathbb H^d$, $d\geq 2$, the first-ordered assumptions are indispensable as well. 
	
	(c) There are more involved, non-Heisenberg-type vector fields which verify also the assumptions of Theorem \ref{fotetel}. Indeed, let $(\mathbb R^5,\mathcal D,g)$ be the sub-Riemannian manifold with the vector fields $X_i,$ $i=1,...,4$ from (\ref{vector-fields}) with $A_1=-2x_2+x_1x_4^2$, $A_2=2x_1$, $A_3=-x_4$, $A_4=x_3+x_1^2x_4.$ In this case we have that the elements from (\ref{A-feltetel}) are $c_{12}=4=-c_{21}$, $c_{34}=2=-c_{43}$, while the rest of $c_{ij}$'s are zero.

	 (d) Note that Theorem \ref{fotetel} can be formulated  on any simply connected open domain instead of the whole $M$. 
\end{remark}

\textit{Organization  of the paper.} In Section \ref{section-2} we prove the Poincar\'e lemma on generic Riemannian manifolds.  As a direct byproduct,  we also state a Saint-Venant lemma on Riemannian manifolds  whose proof is presented in the Appendix (Section \ref{section-appendix}). The Poincar\'e lemma on generic Riemannian manifolds turns to be indispensable in the proof of our main theorem, which will be provided in Section \ref{section-3}. Here, we shall explore  basic properties of the
Riemannian manifolds as the metric compatibility and torsion-freeness (or symmetry) of the Levi-Civita
connection with respect to the Riemannian metric. In fact, we shall reduce our original sub-Riemannian system (defined on the distribution) to a differential system on a Riemannian manifold where we can apply the Riemannian Poincar\'e lemma and Ces\`aro-Volterra integral formula. An elegant computation connects the force vector fields in these two settings, proving in this way relation (\ref{C-V-1}).  In Section \ref{section-examples} we give some examples, the first on the hyperbolic spaces, the second one on Carnot/Heisenberg groups. In Section \ref{remark-section} we formulate some problems for further investigations.  

\section{Poincar\'e lemma on Riemannian manifolds: a local version}\label{section-2}
  Let $(M,g)$ be an $m$-dimensional Riemannian
manifold; here
$(g_{ij})$ are the components of the Riemannian metric $g$ in a given local coordinate system $(x_i)_{i=1,...,m}$. 

Let $u:M\to \mathbb R$ be a $C^1$-functional on $M$; the
differential of $u$ at $x$, denoted by $du(x)$, belongs to the cotangent space $T_x^*M$
and is defined by
\begin{equation}\label{differential}
du(x)(v)=\langle{\bf \nabla}_gu(x),v\rangle_g\ \mbox{for all}\ v\in
T_xM;
\end{equation}
in the sequel, we prefer to use $\langle\cdot,\cdot \rangle_g$ instead of $g$. If the  local components of
$du$ are denoted by $u_k=\partial_{x_k}u$, then the
local components of ${\bf \nabla}_gu$ are $u^i=g^{ik}u_k$; here, $g^{ij}$
are the local components of $g^{-1}=(g_{ij})^{-1}$.

Let  $\Omega\subseteq M$ be an open set and $\textit{\textbf{V}}\in
T\Omega=\cup_{x\in \Omega }T_xM$ be an arbitrary vector field in
$\Omega$ which is represented in local coordinates as
$$\textit{\textbf{V}}=V_k {\partial_{x_k} }.$$ The main result of the present section is the Poincar\'e lemma on Riemannian manifolds.

\begin{theorem}\label{prop-1} Let $(M,g)$ be an $m$-dimensional Riemannian
	manifold and $\Omega\subseteq M$ be a simply connected open set. Given a vector field $\textit{\textbf{V}}\in \textit{\textbf{C}}^1(\Omega, T\Omega)$, the system 
	\begin{equation}\label{gradiens-megoldas}
	{\bf \nabla}_g u=\textit{\textbf{V}}\ {in}\ \Omega
	\end{equation}
	is solvable in $C^2(\Omega)$ if and only if 
  we have
\begin{equation}\label{curl-assumption}
    {\partial_{ x_i}}\tilde V_j={\partial_{ x_j}}\tilde V_i\ \  \ in\ \ \Omega,\ \ {for\ every}\ i,j=1,...,m,
\end{equation}
where $\tilde V_j= g_{jk}V_k.$ 

Moreover, if $x_0\in \Omega$ is fixed and {\rm (\ref{curl-assumption})} holds,  the solution $u:\Omega\to \mathbb R$ for  {\rm (\ref{gradiens-megoldas})} can be obtained by 
\begin{equation}\label{C-V}
u(x)=c_0+\int_0^1 \langle\textit{\textbf{V}}(\gamma(t)),\dot \gamma(t))\rangle_g dt,\ x\in \Omega, 
\end{equation}
where $c_0=u(x_0)\in \mathbb R$ and $\gamma:[0,1]\to \Omega$ is any curve joining $x_0$ with $x.$ 

\end{theorem}

{\it Proof.} "\underline{(\ref{gradiens-megoldas}) implies (\ref{curl-assumption})}" First of all,  
(\ref{gradiens-megoldas}) is equivalent to
$$g^{ik}{\partial_{ x_k}}{u}=V_i,\ \ i=1,...,m.$$ Multiplying both sides by $g_{ji}$, we have that
$${\partial_{ x_j}}{u}=g_{ji}V_i=\tilde V_j,\ j=1,...,m.$$ Deriving these relations,  (\ref{curl-assumption}) yields at once by the symmetry of second-order derivatives.

"\underline{(\ref{curl-assumption}) implies (\ref{gradiens-megoldas})}" We closely follow the proof
from Ciarlet \cite[Theorem 6.17-2]{Ciarlet-2013}. Let $x_0\in
\Omega$ be given and fix $x\in \Omega$. Since $\Omega$ is simply
connected, there exists a path $\gamma:[0,1]\to \Omega$ such that
$\gamma(0)=x_0$ and $\gamma(1)=x.$ If there exists $u\in
C^2(\Omega)$ which satisfies (\ref{gradiens-megoldas}), then the
function $P:[0,1]\to \mathbb R$ defined by $P(t)=u(\gamma(t))$
verifies
$$\frac{dP}{dt}(t)=du(\gamma(t))(\dot \gamma(t))=\langle{\bf \nabla}_gu(\gamma(t)),\dot\gamma(t)\rangle_g,\ t\in [0,1].$$
The latter equation together with the Cauchy data
$P(0)=P_0\in \mathbb R$ provides a unique solution $P:[0,1]\to
\mathbb R$ which depends on the path $\gamma$.

We are going to show that the value $P(1)$ does not depend on
the choice of the path $\gamma$ whenever (\ref{curl-assumption}) holds.  To see this, let $\gamma_0,\gamma_1:[0,1]\to
\Omega$ be two smooth paths such that $\gamma_i(0)=x_0$ and
$\gamma_i(1)=x$, $i\in \{0,1\}.$ Since $\Omega$ is simply connected,
we can find a smooth homotopy $H:[0,1]\times [0,1]\to \Omega$
between $\gamma_0$ and $\gamma_1$, i.e.,
$$H(\cdot,0)=\gamma_0,\ H(\cdot,1)=\gamma_1,\ $$
$$H(0,\lambda)=x_0,\ H(1,\lambda)=x,\ \forall \lambda\in [0,1]. $$
For every $\lambda\in [0,1]$, let $P(\cdot,\lambda):[0,1]\to \mathbb
R$ be the unique solution of the Cauchy problem
\[ \displaystyle \   \left\{ \begin{array}{lll}
 \frac{\partial P}{\partial t}(t,\lambda)=\left\langle \textit{\textbf{V}}(H(t,\lambda)),\frac{\partial H}{\partial t}(t,\lambda)\right\rangle_g,\  &\mbox{for} & t\in [0,1]; \\
P(0,\lambda)=P_0\in \mathbb R.
 \end{array}\right. \eqno{({\mathcal C}_{\lambda})}\]
 We claim that
 \begin{equation}\label{claim}
\frac{\partial P}{\partial \lambda}(1,\lambda)=0\ \ {\rm for\ every}\
\lambda\in [0,1].
 \end{equation}
To see this, let us consider the function $\sigma:[0,1]\times
[0,1]\to \mathbb R$ defined by
$$\sigma(t,\lambda)=\frac{\partial P}{\partial \lambda}(t,\lambda)-\left\langle \textit{\textbf{V}}(H(t,\lambda)),\frac{\partial H}{\partial \lambda}(t,\lambda)\right\rangle_g.$$
Since the Levi-Civita connection is compatible with the  Riemannian
metric, it follows by do Carmo \cite[Proposition 3.2]{doCarmo} that
$$\frac{\partial\sigma}{\partial t}(t,\lambda)=\frac{\partial}{\partial t}\left(\frac{\partial P}{\partial \lambda}\right)(t,\lambda)-
\left\langle \frac{D\textit{\textbf{V}}}{\partial t}(H(t,\lambda)),\frac{\partial
H}{\partial \lambda}(t,\lambda)\right\rangle_g- \left\langle
\textit{\textbf{V}}(H(t,\lambda)),\frac{D}{\partial t}\frac{\partial H}{\partial
\lambda}(t,\lambda)\right\rangle_g,$$ where $D$ denotes the
covariant derivation on $(M,g)$.  Concerning the latter term, we know from the torsion-freeness of the Levi-Civita connection on $(M,g)$ that 
\begin{equation}\label{symmetry-1}
\frac{D}{\partial t}\frac{\partial H}{\partial
\lambda}(t,\lambda)=\frac{D}{\partial \lambda}\frac{\partial
H}{\partial t}(t,\lambda),
\end{equation}
see do
Carmo \cite[Lemma 3.4]{doCarmo}. The sophisticated part is to show that
\begin{equation}\label{symmetry-2}\left\langle
\frac{D\textit{\textbf{V}}}{\partial t}(H(t,\lambda)),\frac{\partial H}{\partial
\lambda}(t,\lambda)\right\rangle_g=\left\langle \frac{D\textit{\textbf{V}}}{\partial
\lambda}(H(t,\lambda)),\frac{\partial H}{\partial
t}(t,\lambda)\right\rangle_g.
\end{equation}
To prove (\ref{symmetry-2}) we recall the following well known
facts: if $\textit{\textbf{W}}=(w_1,...,w_m)$ is a vector field along a path $(x)$,
its covariant derivative can be expressed by
$$\frac{D\textit{\textbf{W}}}{dt}=\left(\frac{dw_k}{dt}+\Gamma_{ij}^kw_j\frac{dx_i}{dt}\right)\partial _{x_k},$$
where $\Gamma_{ij}^k$ are the Christofel symbols for which we have
\begin{equation}\label{Christo}
g_{ks}\Gamma_{ij}^k=\frac{1}{2}\left({\partial_{x_i} g_{js}}+\partial_{x_j}{ g_{is}}-\partial_{x_s}
g_{ij}\right).
\end{equation}
Coming back to (\ref{symmetry-2}), we have
\begin{eqnarray*}
  LHS &:=& \left\langle \frac{D\textit{\textbf{V}}}{\partial t}(H(t,\lambda)),\frac{\partial
H}{\partial
\lambda}(t,\lambda)\right\rangle_g=g_{kj}\left({\partial_{x_i}}{
V_k}\frac{\partial H_i}{\partial
t}+\Gamma_{il}^kV_l\frac{\partial H_i}{\partial
t}\right)\frac{\partial H_j}{\partial \lambda} \\
   &=&g_{kj}\left({\partial_{x_i}}{
   	V_k}+\Gamma_{il}^kV_l\right)\frac{\partial
H_i}{\partial t}\frac{\partial H_j}{\partial \lambda}.
\end{eqnarray*}
In a similar way,
\begin{eqnarray*}
  RHS &:=& \left\langle \frac{D\textit{\textbf{V}}}{\partial \lambda}(H(t,\lambda)),\frac{\partial
H}{\partial t}(t,\lambda)\right\rangle_g=g_{ki}\left({\partial_{x_j}}{
V_k}\frac{\partial H_j}{\partial
\lambda}+\Gamma_{jl}^kV_l\frac{\partial H_j}{\partial
\lambda}\right)\frac{\partial H_i}{\partial t} \\
   &=&g_{ki}\left({\partial_{x_j}}{
   	V_k}+\Gamma_{jl}^kV_l\right)\frac{\partial
H_i}{\partial t}\frac{\partial H_j}{\partial \lambda}.
\end{eqnarray*}
Therefore, we have that
\begin{eqnarray*}
  (\ref{symmetry-2})\ {\rm holds} &\Leftrightarrow&  LHS-RHS=0\\
   &\Leftrightarrow& \left[g_{kj}\left({\partial_{x_i}}{
   	V_k}+\Gamma_{il}^kV_l\right)-g_{ki}\left({\partial_{x_j}}{
   	V_k}+\Gamma_{jl}^kV_l\right)\right]\frac{\partial
H_i}{\partial t}\frac{\partial H_j}{\partial \lambda}=0 \\
   &\Leftrightarrow& \left[g_{kj}{\partial_{x_i}}{
   	V_k}-g_{ki}{\partial_{x_j}}{
   	V_k}+(g_{kj}\Gamma_{il}^k-g_{ki}\Gamma_{jl}^k)V_l\right]\frac{\partial
H_i}{\partial t}\frac{\partial H_j}{\partial \lambda}=0 \\
   &\stackrel{(\ref{Christo})}{\Leftrightarrow}&\left[g_{kj}{\partial_{x_i}}{
   	V_k}-g_{ki}{\partial_{x_j}}{
   	V_k}+\left(\partial_{x_i} g_{lj}-\partial_{x_j}
g_{li}\right)V_l\right]\frac{\partial H_i}{\partial
t}\frac{\partial H_j}{\partial \lambda}=0\\
 &{\Leftrightarrow}&\left[g_{kj}{\partial_{x_i}}{
 	V_k}-g_{ki}{\partial_{x_j}}{
 	V_k}+\left(\partial_{x_i} g_{kj}-\partial_{x_ij}
g_{ki}\right)V_k\right]\frac{\partial H_i}{\partial
t}\frac{\partial H_j}{\partial \lambda}=0\\
&{\Leftrightarrow}&\left[{\partial_{x_i}}(g_{jk}V_k)-\partial_{x_j}(g_{ik}V_k)\right]\frac{\partial H_i}{\partial t}\frac{\partial
H_j}{\partial \lambda}=0,
\end{eqnarray*}
where the latter relation holds true due to (\ref{curl-assumption}).
Consequently, by relations (\ref{symmetry-1}), (\ref{symmetry-2})
and the Cauchy problem $(\mathcal C_\lambda)$ we have
\begin{eqnarray*}
  \frac{\partial\sigma}{\partial t}(t,\lambda) &=& \frac{\partial}{\partial \lambda}\left(\frac{\partial P}{\partial t}\right)(t,\lambda)-
\left\langle \frac{D\textit{\textbf{V}}}{\partial
\lambda}(H(t,\lambda)),\frac{\partial H}{\partial
t}(t,\lambda)\right\rangle_g- \left\langle
\textit{\textbf{V}}(H(t,\lambda)),\frac{D}{\partial \lambda}\frac{\partial H}{\partial
t}(t,\lambda)\right\rangle_g \\
   &=& \frac{\partial}{\partial \lambda}\left(\frac{\partial P}{\partial t}(t,\lambda)-\left\langle \textit{\textbf{V}}(H(t,\lambda)),\frac{\partial H}{\partial t}(t,\lambda)\right\rangle_g\right) \\
   &=&0,
\end{eqnarray*}
i.e., $t\mapsto \sigma(t,\lambda)$ is constant. Since $P(0,\lambda)=P_0\in \mathbb R$ and $H(0,\lambda)=x_0$, it turns out that 
 $$\sigma(0,\lambda)=\frac{\partial P}{\partial \lambda}(0,\lambda)-\left\langle \textit{\textbf{V}}(H(0,\lambda)),\frac{\partial H}{\partial \lambda}(0,\lambda)\right\rangle_g=0\  \ {\rm for\ every}\ \lambda\in [0,1].$$ In particular, $$0=\sigma(1,\lambda)=\frac{\partial P}{\partial \lambda}(1,\lambda)-\left\langle \textit{\textbf{V}}(H(1,\lambda)),\frac{\partial H}{\partial \lambda}(1,\lambda)\right\rangle_g.$$ Since $H(1,\lambda)=x_0$ for every $\lambda\in [0,1]$, it follows the claim
(\ref{claim}), showing that the value $P(1)$ is not depending on the
particular choice of the path.

For every $x\in \Omega,$ let $u:\Omega\to \mathbb R$ be defined by
$$u(x)=P(1),$$
where $P$ is the unique solution to the Cauchy problem $({\mathcal C}_{\lambda})$ having the
initial data $P(0)=P_0$ and using any path joining $x_0$ and $x$;
thus, the function $u$ is well defined.

 To conclude the proof,
we show the validity of (\ref{gradiens-megoldas}). Let $x\in
\Omega$ and $v\in T_xM$ be arbitrarily fixed elements. Let
$\gamma:[0,1]\to \Omega$ be a path such that $\gamma(0)=x_0$,
$\gamma(1)=x$ and $\dot \gamma(1)=v\in T_xM,$ and let $P$ be the solution of
the Cauchy problem associated to this path, thus,
$P(t)=u(\gamma(t))$. Therefore, the latter relation yields that
$$\frac{dP}{dt}(t)=\langle{\bf \nabla}_gu(\gamma(t)),\dot\gamma(t)\rangle_g,\ t\in [0,1].$$
On the other hand, by the Cauchy problem we have
$$\frac{dP}{dt}(t)=\left\langle \textit{\textbf{V}}(\gamma(t)),\dot\gamma(t)\right\rangle_g,\ t\in [0,1].$$
Accordingly, for the moment $t=1$, it follows that
$$\left\langle \nabla_gu(x),v\right\rangle_g=\left\langle
\textit{\textbf{V}}(x),v\right\rangle_g$$ and the arbitrariness of $v\in T_xM$
concludes the proof of (\ref{gradiens-megoldas}).

If $\gamma:[0,1]\to \Omega$ is any path joining the points $x_0$ and $x$, the Ces\`aro-Volterra path integral formula easily follows as
$$u(x)-u(x_0)=\int_0^1 \frac{d}{dt}u(\gamma(t))dt=\int_0^1\langle{\bf \nabla}_gu(\gamma(t)),\dot\gamma(t)\rangle_gdt=\int_0^1\langle\textit{\textbf{V}}(\gamma(t)),\dot\gamma(t)\rangle_gdt,$$
which is precisely (\ref{C-V}). \hfill $\square$

\begin{remark}\rm
	Poincar\'e's lemma can be also proved by using 1-forms, see e.g. Abraham, Marsden and Ratiu \cite{AMR}. However, we preferred here a direct proof  based on local coordinates for two reasons: (a) it highlights the importance of the Riemannian structure, i.e., the metric compatibility and torsion-freeness  of the Levi-Civita	connection, which is not valid anymore on non-Riemannian Finsler settings (see  Section \ref{remark-section} for details); (b) the proof provides directly  a Ces\`aro-Volterra path integral formula. 
\end{remark}

As a byproduct of the Poincar\'e lemma (Theorem \ref{prop-1}), we state a Saint-Venant lemma on generic Riemannian manifolds;  its proof is sketched in the Appendix. To present it, fix $\textit{\textbf{e}}_i\in T\Omega$, $i=1,...,m$, and assume that they can be
represented as
$$\textit{\textbf{e}}_i=e_{ik}{\partial_{x_k}}.$$
The $m$-vector field $\textit{\textbf{e}}=(\textit{\textbf{e}}_1,...,\textit{\textbf{e}}_m)\in \textit{\textbf{C}}^2(\Omega,T\Omega^m)$ is called {\it symmetric}
if $e_{ij}=e_{ji}\in C^2(\Omega)$ for every $i,j=1,...,m.$

\begin{proposition}\label{Saint-Venant} Let $(M,g)$ be an $m$-dimensional Riemannian
	manifold and $\Omega\subseteq M$ be a simply connected open set.
	Given  $\textit{\textbf{e}}=(\textit{\textbf{e}}_1,...,\textit{\textbf{e}}_m)\in \textbf{C}^2(\Omega,T\Omega^m)$  a
	symmetric $m$-vector field on $\Omega$, 
	the system 
		\begin{equation}\label{szim-deriv}
	\nabla_{s,g} \textbf{V}=\textit{\textbf{e}}\ \ {in}\ \ \Omega,
	\end{equation}
	has a vector field solution $\textbf{V}=(V_1,...,V_m)\in \textbf{C}^3(\Omega,\mathbb R^m),$ 
	where the components of the symmetric gradient $\nabla_{s,g} \textbf{V}$  are given by
	$$\frac{1}{2}\left(\partial_{x_i}(g_{jk}V_k)+\partial_{x_j}(g_{ik}V_k)\right),\ \
	i,j=1,...,m,$$ if and only if the Saint-Venant
	compatibility relations hold $($in local coordinate system$)$ in
	$\Omega$, i.e.,
	\begin{equation}\label{SV-relation}
	\partial^2_{x_lx_j}e_{ik}+\partial^2_{x_kx_i}e_{jl}-\partial^2_{x_lx_i}e_{jk}-\partial^2_{x_jx_k}e_{il}=0,\ \ i,j,k,l=1,...,m.
	\end{equation}

 Moreover, if $x_0\in M$ is fixed and {\rm (\ref{SV-relation})} holds, then the solution of {\rm(\ref{szim-deriv})} is obtained by  $$V_k=g^{ks}u_s,\ \ k=1,...,m,$$
 where $$u_i(x)=c_0^i+\int_0^1\langle \textit{\textbf{U}}_i(\gamma(t)),\dot\gamma(t)\rangle_gdt,\ \ x\in \Omega,$$
 with $\textit{\textbf{U}}_i=g^{ls}(p_{is}+e_{is}){\partial}_{x_l},$ 
 $$p_{ij}(x)=c_0^{ij}+\int_0^1\langle \textit{\textbf{W}}_{ij}(\gamma(t)),\dot\gamma(t)\rangle_gdt,\ \ x\in \Omega,$$ and 
 $\textit{\textbf{W}}_{ij}=g^{ls}\left({\partial_{ x_j}e_{is}}{}-{\partial_{ x_i}e_{js}}\right){\partial}_{x_l},$  for some numbers $c_0^s,$ $c_0^{ij}$ and the curve $\gamma:[0,1]\to \Omega$ is arbitrary fixed joining $x_0$ with $x\in \Omega.$ 
\end{proposition}

\begin{remark}\rm
(a)	Note that 	$\nabla_{s,g} \textit{\textbf{\textbf{V}}}$  is a kind of {\it symmetric Lie derivative}
	of the  vector field $\textit{\textbf{\textbf{V}}}$ with respect to the Riemannian metric $g$; indeed, the latter notion appears in   Chen
	and Jost \cite[p. 518]{Chen-Jost}, where  $\nabla_{s,g} \textit{\textbf{\textbf{V}}}$ is an
	$\mathcal L-$type tensor of the form 
	$$ \nabla_{s,g} \textit{\textbf{\textbf{V}}}=\frac{1}{2}\left(g_{jk}{\partial_{x_i} V_k}+g_{ik}{\partial_{x_j} V_k}+C_{ijk}V_k\right)dx_i\otimes dx_j.$$ In our setting, the elements $C_{ijk}$ are 
	expressed by means of the Christoffel symbols as
	$$C_{ijk}={\partial_{x_i} g_{jk}}+{\partial_{x_j}} g_{ik}=g_{lj}\Gamma_{ki}^l+g_{li}\Gamma_{kj}^l+2g_{lk}\Gamma_{ij}^l.$$	
	
	(b) Proposition \ref{Saint-Venant} provides a curved version of the Saint-Venant lemma; further curvilinear versions of the Saint-Venant lemma can be found in the papers by Ciarlet,  Gratie, Mardare and  Shen \cite{Ciarlet-Mardare-2}, Ciarlet and  Mardare \cite{Ciarlet-Mardare-0}, and Ciarlet, Mardare and Shen \cite{Ciarlet-Mardare}.
\end{remark}


\section{Proof of Theorem \ref{fotetel}}\label{section-3}

In order to prove Theorem \ref{fotetel}, we first recall some basic notions from the theory of sub-Riemannian manifolds; for further details, see Agrachev,  Barilari and Boscain \cite{ABB}, Calin and Chang \cite{Calin-book} and Figalli and Rifford \cite{FR}.  

Let $M$ be a smooth connected $(n+1)$-dimensional manifold $(n\geq 2)$, $\mathcal D$ be a smooth nonholonomic  distribution of rank $m\leq n$ on $M$ (i.e., a rank $m$ subbundle of the tangent bundle $TM$) and $g$ be a Riemannian metric on $\mathcal D$. Without loss of generality, we may assume that $g$ is defined on the whole tangent bundle $TM$  (not necessarily in a unique way); we shall keep the same notation of $g$ on $TM$.  The triplet $(M,\mathcal D,g)$ is a sub-Riemannian manifold. As usual, the distribution $\mathcal D$ is said to be nonholonomic if for every $x\in M$ there exists an $m$-tuple $X_1^x,...,X_m^x$ of smooth vector fields on a neighborhood $N_x$ of $x$ such that all the Lie brackets generated by these vectors at $y$ generate $T_yM$  for every $y\in N_x.$ A curve $\gamma:[0,1]\to M$ is horizontal with respect to $\mathcal D$ if it belongs to $W^{1,2}([0,1];M)$ and $\dot \gamma(t)\in \mathcal D(\gamma(t))$ for a.e. $t\in [0,1]$. If $\mathcal D$ is nonholonomic on $M,$ by the Chow-Rashewsky theorem, every two points of $M$ can be joined by a horizontal path. Let $\Gamma(\mathcal D)$ be the set of horizontal vector fields on $M,$ and $\mathcal F(M)$ be the set of smooth functions on $M.$ If $u\in \mathcal F(M)$, the horizontal gradient $\nabla_Hu\in \Gamma(\mathcal D)$ of $u$  is defined by $g(\nabla_Hu,X)=X(u)$ for every  $X\in \Gamma(\mathcal D).$  

Now, let us put ourselves into the context of Theorem \ref{fotetel}. Accordingly, let $(M,\mathcal D,g)$ be an $(n+1)$-dimensional sub-Riemannian manifold $(n\geq 2)$,  and the  rank $n$ distribution $\mathcal D$ in a local coordinate system $(x_i)_{i=1,...,n+1}$ formed by the vector fields given in (\ref{vector-fields}) and   
verifying 
 (\ref{A-feltetel}).  Since 
\begin{eqnarray*}
X_iX_j&=&(\partial_{x_i}+A_i\partial_{ x_{n+1}})(\partial_{x_j}+A_j\partial_{ x_{n+1}})\\&=&\partial^2_{x_ix_j}+\partial_{x_i}A_j\partial_{ x_{n+1}}+A_j\partial^2_{ x_ix_{n+1}}+A_i\partial^2_{x_j x_{n+1}}+A_iA_j\partial^2_{x_{n+1}},
\end{eqnarray*} 
   by (\ref{A-feltetel})  we obtain (\ref{Xi-k-c-ik}), i.e., 
$$
[X_i,X_j]=X_iX_j-X_jX_i=(\partial_{x_i}A_j-\partial_{x_i}A_j)\partial_{ x_{n+1}}=c_{ij}\partial_{x_{n+1}}\ \ {\rm for\ every}\  i,j=1,...,n.
$$
Therefore, since  
$I_0=\{(i,j): c_{ij}\neq 0\}\neq \emptyset,$ the 
distribution $\mathcal D$ is nonholonomic on $M$.

Let $\textit{\textbf{a}}\in \Gamma(\mathcal D)$ be fixed. The system (\ref{sub-riem-equation}), i.e., $${\bf\nabla}_H u=\textit{\textbf{a}},$$ in local coordinates reads as
\begin{equation}\label{lokal-szub-riemann}
X_j(u)=g_{ij}a_i=:\tilde a_j,\ \ j=1,...,n,
\end{equation}
where $g_{ij}=g(X_i,X_j)$ and $\textit{\textbf{a}}=a_iX_i.$
With this preparatory part in our mind, we now present the\\

{\it Proof of Theorem \ref{fotetel}}. "\underline{(\ref{sub-riem-equation}) implies (\ref{full-jellemzes})\&(\ref{full-jellemzes-2})}" Assume that the sub-Riemannian system (\ref{sub-riem-equation}) has a solution $u\in \mathcal F(M)$. First, by (\ref{Xi-k-c-ik}) applied to $u$, we have that 
$$[X_i,X_j]u=c_{ij}\partial_{x_{n+1}}u,\ \ i,j=1,...,n.$$ This relation and  (\ref{lokal-szub-riemann}) give that
\begin{equation}\label{kell-elso}
X_i \tilde a_j-X_j \tilde a_i=c_{ij}\partial_{x_{n+1}}u,\ \ i,j=1,...,n.
\end{equation}

If $\partial_{x_{n+1}}u(x)=0$ for some $x\in M$, then $X_i \tilde a_j(x)-X_j \tilde a_i(x)=0$ for every $i,j=1,...,n,$ thus (\ref{full-jellemzes}) clearly holds. If  $\partial_{x_{n+1}}u(x)\neq 0$ for some $x\in M$, then by writing the  relation (\ref{kell-elso}) for $(k,l)$ instead of $(i,j)$, and eliminating $\partial_{x_{n+1}}u(x)\neq 0$, we obtain (\ref{full-jellemzes}). 

 Deriving (\ref{kell-elso}) with respect to the vector field $X_k$, $k=1,...,n$, and taking into account that $[X_k,\partial_{x_{n+1}}]=X_k\partial_{x_{n+1}}-\partial_{x_{n+1}}X_k=0$, it turns out by (\ref{lokal-szub-riemann}) and (\ref{Xi-k-c-ik}) that 
 $$X_kX_i \tilde a_j-X_kX_j \tilde a_i=c_{ij}X_k\partial_{x_{n+1}}u=c_{ij}\partial_{x_{n+1}}X_ku=[X_i,X_j]\tilde a_k,$$
which is precisely relation (\ref{full-jellemzes-2}). 

"\underline{(\ref{full-jellemzes})\&(\ref{full-jellemzes-2}) imply (\ref{sub-riem-equation})}" 
Since $I_0\neq \emptyset$, let $(i_0,j_0)\in I_0$ and introduce the function  $$\tilde a=\frac{X_{i_0}\tilde a_{j_0}-X_{j_0}\tilde a_{i_0}}{c_{i_0j_0}},$$
where $\tilde a_j=g_{ij}a_i$. With these notations, we consider the system
\begin{equation}\label{S-system}
	\left\{
\begin{array}{lll}
\partial_{x_j}u=\tilde a_j-A_j\tilde a&
{\rm for} & j=1,...,n;\\
 \partial_{x_{n+1}}u=\tilde a.
\end{array}
\right.
\end{equation}
Let $$\tilde V_j=\tilde a_j-A_j\tilde a\ (j=1,...,n)\ \ {\rm and}\ \  \tilde V_{n+1}=\tilde a;$$
we are going to prove that 
\begin{equation}\label{megnezni}
\partial_{x_i}\tilde V_j=\partial_{x_j}\tilde V_i, \ \ i,j=1,...,n+1.
\end{equation}
To do this, we distinguish three cases: 

\underline{Case 1}: $i=j=n+1$;  (\ref{megnezni}) trivially holds. 

\underline{Case 2}: $i\in \{1,...,n\}$ and $j=n+1.$ On one hand, (\ref{megnezni}) is equivalent to $\partial_{x_i}\tilde a=\partial_{x_{n+1}}(\tilde a_i-A_i\tilde a)$, which can be written as $X_i\tilde a=\partial_{x_{n+1}}\tilde a_i$. On the other hand, by the definition of $\tilde a$,  (\ref{full-jellemzes-2}) and (\ref{Xi-k-c-ik}) we have that
$$X_i\tilde a=\frac{X_iX_{i_0}\tilde a_{j_0}-X_iX_{j_0}\tilde a_{i_0}}{c_{i_0j_0}}=\frac{[X_{i_0},X_{j_0}]\tilde a_i}{c_{i_0j_0}}=\partial_{x_{n+1}}\tilde a_i,$$ which is the required relation. 

\underline{Case 3}: $i,j\in \{1,...,n\}$. We have the following chain of equivalences: 
\begin{eqnarray*}
(\ref{megnezni})\ {\rm holds}&\Leftrightarrow& \partial_{x_i}\tilde a_j-\tilde a \partial_{x_i}A_j-A_j\partial_{x_i}\tilde a=\partial_{x_j}\tilde a_i-\tilde a \partial_{x_j}A_i-A_i\partial_{x_j}\tilde a\\
&\Leftrightarrow& \partial_{x_i}\tilde a_j-\tilde a \partial_{x_i}A_j-A_jX_i\tilde a=\partial_{x_j}\tilde a_i-\tilde a \partial_{x_j}A_i-A_iX_j\tilde a\\
&\Leftrightarrow& \partial_{x_i}\tilde a_j-A_jX_i\tilde a=\partial_{x_j}\tilde a_i+ c_{ij}\tilde a -A_iX_j\tilde a\\
&\Leftrightarrow& \partial_{x_i}\tilde a_j-A_j\frac{[X_{i_0},X_{j_0}]\tilde a_i}{c_{i_0j_0}}=\partial_{x_j}\tilde a_i+ c_{ij}\tilde a -A_i\frac{[X_{i_0},X_{j_0}]\tilde a_j}{c_{i_0j_0}} \ \ \ \ \ ({\rm by}\ (\ref{full-jellemzes-2})) \\
&\Leftrightarrow& \partial_{x_i}\tilde a_j-A_j\partial_{x_{n+1}}\tilde a_i=\partial_{x_j}\tilde a_i+ c_{ij}\tilde a -A_i\partial_{x_{n+1}}\tilde a_j \ \ \ \ \ \ \ \ \ \ \ \ \ \ \ \  \ ({\rm by}\ (\ref{Xi-k-c-ik})) \\
&\Leftrightarrow& X_i\tilde a_j-X_j\tilde a_i=c_{ij}\tilde a.
\end{eqnarray*}
 By the definition of $\tilde a$, let us observe that the latter relation is nothing but (\ref{full-jellemzes}) with the choice $(k,l)=(i_0,j_0)$, which concludes the proof of (\ref{megnezni}).
   
 According to Theorem \ref{prop-1} (applied for $(M,\tilde g)$ with  $\tilde g_{ij}=g(\partial_{ x_i},\partial_{ x_j})$, $i,j=1,...,n+1$) and relation (\ref{megnezni}), it turns out that the system (\ref{S-system}) has a solution in $C^2(M)$, which can be obtained by 
 \begin{equation}\label{egyenlet-riemann}
  u(x)=c_0+\int_0^1 \langle\textit{\textbf{V}}(\gamma(t)),\dot \gamma(t))\rangle_{\tilde g} dt,\ x\in M,
 \end{equation}
 where  $\textit{\textbf{V}}=\sum_{i=1}^{n+1}V_i\partial_{ x_i}$ with $V_i=\sum_{j=1}^{n+1}\tilde g^{ij}\tilde V_j$ and ${\tilde g}^{ij}=({\tilde g}_{ij})^{-1}$;  here, $\gamma:[0,1]\to M$ is any curve joining an $x_0\in M$ with $x\in M$, with $c_0=u(x_0)$. 
   
   By (\ref{S-system}) we clearly have for every $j=1,...,n$ that 
   $$X_j(u)=\partial_{x_j}u+A_j\partial_{x_{n+1}}u=(\tilde a_j-A_j\tilde a)+A_j\tilde a=\tilde a_j,$$
   which is equivalent to ${\bf\nabla}_H u=\textit{\textbf{a}}$, see (\ref{lokal-szub-riemann}), i.e., $u\in C^2(M)$ is a solution to  (\ref{sub-riem-equation}). 
   
   It remains to prove the sub-Riemannian Ces\`aro-Volterra path integral formula (\ref{C-V-1}). To do this, let us fix an arbitrary horizontal path $\gamma:[0,1]\to M$, joining $x_0$ with $x\in M$. If $\gamma$  has the local representation $\gamma=(\gamma_1,...,\gamma_{n+1})$, its horizontality means that $$\dot\gamma_{n+1}=\sum_{k=1}^n A_k\dot\gamma_k.$$ Considering every term at the moment $t\in [0,1]$ in the following computations, we have  
   \begin{eqnarray*}
   \langle\textit{\textbf{V}}(\gamma(t)),\dot \gamma(t))\rangle_{\tilde g}&=&\sum_{i,k=1}^{n+1}\tilde g_{ik}V_i\dot \gamma_k=\sum_{i,k,j=1}^{n+1}\tilde g_{ik}\tilde g^{ij}\tilde V_j\dot \gamma_k=\sum_{k=1}^{n+1}\left(\sum_{j=1}^{n+1}\left(\sum_{i=1}^{n+1}\tilde g_{ik}\tilde g^{ij}\right)\tilde V_j\right)\dot{\gamma_k}\\&=&\sum_{k=1}^{n+1}\left(\sum_{j=1}^{n+1}\delta_{kj}\tilde V_j\right)\dot{\gamma_k}=\sum_{k=1}^{n+1}\tilde V_k\dot{\gamma_k}=\sum_{k=1}^{n}\tilde V_k\dot{\gamma_k}+\tilde V_{n+1}{\dot\gamma_{n+1}}\\&=&\sum_{k=1}^{n}(\tilde V_k+A_k\tilde V_{n+1})\dot\gamma_k\\&=&
   \sum_{k=1}^{n}(\tilde a_k-A_k\tilde a+A_k\tilde a)\dot\gamma_k=\sum_{k=1}^{n}\tilde a_k\dot\gamma_k =\sum_{k=1}^{n}g_{ik} a_i\dot\gamma_k\\&=&
  g(\textit{\textbf{a}}(\gamma(t)),\dot \gamma(t)).
   \end{eqnarray*}
   Thus, by (\ref{egyenlet-riemann}) and the latter computation we obtain  (\ref{C-V-1}), which concludes our proof. 
\hfill $\square$

\section{Examples}\label{section-examples}

In this section we provide some computational examples as applications to Theorems \ref{fotetel} \& \ref{prop-1} and Proposition \ref{Saint-Venant}, respectively.

\subsection{Hyperbolic space} Let $\mathbb B^m=\{x\in \mathbb R^m:|x|<1\}$ be the set 
endowed with the Riemannian metric $$g_{\rm
	hyp}(x)=(g_{ij}(x))_{i,j={1,...,m}}=p(x)^2\delta_{ij},$$ where
$$p(x)=\frac{2}{1-|x|^2}.$$ The pair $(\mathbb
B^m,g_{\rm hyp})$ is a model of the $m$-dimensional hyperbolic space with constant
sectional curvature $-1$.

\begin{example}\rm 
We solve the problem 
	\begin{equation}\label{hyper}
	\nabla_{g_{\rm
			hyp}}u=\frac{x}{p}\ {\rm in}\ \mathbb B^m,
	\end{equation}where $\nabla_{g_{\rm
			hyp}}$ denotes the hyperbolic gradient. 
	\end{example}
A direct computation shows  that $\partial_{x_i}(px_j)=\partial_{x_j}(px_i)$ for every $i,j=1,...,m,$ thus we may apply Theorem \ref{prop-1} on $(\mathbb
B^m,g_{\rm hyp})$, which implies the solvability of (\ref{hyper}).  
Moreover, if $\gamma:[0,1]\to \mathbb
B^m$ is $\gamma(t)=tx$ with an arbitrarily fixed $x\in \mathbb
B^m$, the solution $u$ can be obtained as
\begin{eqnarray*}
u(x)&=&c_0+\int_0^1\displaystyle\big\langle\frac{\gamma(t)}{p(\gamma(t))},\dot\gamma(t)\big\rangle_{g_{\rm hyp}} dt=c_0+\int_0^1{p(\gamma(t))}\langle{\gamma(t)},\dot\gamma(t)\rangle dt\\&=&c_0+2\int_0^1\frac{|x|^2t}{1-|x|^2t^2} dt=c_0
-\ln(1-|x|^2)\\&=&c_0
+\ln (p(x)/2),
\end{eqnarray*}
for any $c_0\in \mathbb R.$\\

For simplicity, in the next example we consider only the hyperbolic plane $(\mathbb
B^2,g_{\rm hyp})$. 
\begin{example}\rm 
We solve the problem 
	\begin{equation}\label{hyper-2}
	\nabla_{s,g_{\rm
			hyp}}\textit{\textbf{V}}=\textit{\textbf{e}}\ {\rm on}\ \mathbb B^2,
	\end{equation}where $\nabla_{s,g_{\rm
			hyp}}$ denotes the symmetric hyperbolic  gradient and $\textit{\textbf{e}}=(\textit{\textbf{e}}_1,\textit{\textbf{e}}_2)\in \textit{\textbf{C}}^\infty(\mathbb B^2,(T\mathbb B^2)^2)$ has the components $\textit{\textbf{e}}_1=-\frac{x_1}{p}{\partial_{x_2}}$ and $\textit{\textbf{e}}_2=-\frac{1}{p}(x_1{\partial_{x_1}}+2x_2{\partial_{x_2}}).$
\end{example}

First, we have $e_{11}=0,$ $e_{12}=e_{21}=-\frac{x_1}{p}$ and $e_{22}=-\frac{2x_2}{p}.$ It is easily seen that the Saint-Venant relations (\ref{SV-relation}) are verified; for instance, if $i=k=1$ and $j=l=2$  then the components in (\ref{SV-relation}) are $\partial_{x_2x_2}^2 e_{11}=0$, $\partial_{x_1x_1}^2 e_{22}=2x_2$ and $\partial_{x_1x_2}^2 e_{12}=x_2$. Therefore, we may apply Proposition \ref{Saint-Venant}, guaranteeing the solvability of (\ref{hyper-2}). By keeping the same notations as in Proposition \ref{Saint-Venant}, since $g_{\rm
	hyp}^{-1}=p(x)^{-2}\delta_{ij}$,  after some computation it turns out that 
$$\textit{\textbf{W}}_{11}=\textit{\textbf{W}}_{22}=0\ \ {\rm and}\ \  \textit{\textbf{W}}_{12}=-\textit{\textbf{W}}_{21}=\frac{1}{2p^2}(1-|x^2|-2x_1^2)\partial_{x_1}-\frac{x_1x_2}{p^2}\partial_{x_2}.$$
Accordingly, for every $x\in \mathbb B^2$  on has $p_{11}(x)=c_0^{11},$ $p_{22}(x)=c_0^{22}$ for some $c_0^{11},c_0^{22}\in \mathbb R$ and if we fix $\gamma:[0,1]\to \mathbb B^2$ with $\gamma(t)=tx=(tx_1,tx_2)$, then 
$$p_{12}(x)=-p_{21}(x)=c_0^{12}+\int_0^1\langle \textit{\textbf{W}}_{12}(\gamma(t)),\dot\gamma(t)\rangle_{g_{\rm
		hyp}}dt=c_0^{12}+\frac{1}{2}(x_1-x_1^3-x_1x_2^2),$$ for some $c_0^{12}\in \mathbb R.$
Thus, 
$$\textit{\textbf{U}}_1=\frac{1}{p^2}(c_0^{11}\partial_{x_1}+c_0^{12}\partial_{x_2}),$$ and 
$$\textit{\textbf{U}}_2=\frac{1}{p^2}\left((-c_0^{12}-x_1+x_1^3+x_1x_2^2)\partial_{x_1}+(c_0^{22}-x_2+x_1^2x_2+x_2^3)\partial_{x_2}\right).$$
Therefore, for every $x\in \mathbb B^2$, if $\gamma:[0,1]\to \mathbb B^2$ is again the curve given by $\gamma(t)=tx=(tx_1,tx_2)$, then the latter vector fields provide the functions 
$$u_1(x)=c_0^1+\int_0^1\langle \textit{\textbf{U}}_{1}(\gamma(t)),\dot\gamma(t)\rangle_{g_{\rm
		hyp}}dt=c_0^1+c_0^{11}{x_1}+c_0^{12}{x_2},$$
	and 
$$u_2(x)=c_0^2+\int_0^1\langle \textit{\textbf{U}}_{2}(\gamma(t)),\dot\gamma(t)\rangle_{g_{\rm
		hyp}}dt=c_0^2-\frac{1}{4}-c_0^{12}{x_1}+c_0^{22}{x_2}+\frac{1}{p^2(x)}.$$
Consequently, $\textit{\textbf{V}}=(V_1,V_2)$ is a solution of (\ref{hyper-2}), where 
$V_i=\frac{u_i}{p^2}$, $i=1,2,$ with $c_0^{11}=c_0^{22}=0$ and $c_0^1,$ $c_0^2$ and $c_0^{12}$ arbitrarily fixed.

\subsection{Carnot and Heisenberg groups} Let $\mathbb G$ be an $(n+1)$-dimensional corank 1 Carnot group with the Lie algebra $\mathfrak g=\mathfrak g_1\oplus \mathfrak g_2$, where  dim$\mathfrak g_1=n$ and dim$\mathfrak g_2=1$. Usually, the operation on $\mathfrak g$ (in exponential coordinates on $\mathbb R^n\times \mathbb R$) is given by 
$$x\circ y=\left(x_1+y_1,...,x_{n}+y_{n}, x_{n+1}+y_{n+1}-\frac{1}{2}\sum_{i,j=1}^n \mathcal A_{ij}x_jy_i\right ),$$
where $x=(x_1,...,x_{n+1})$, $y=(y_1,...,y_{n+1})$, and  
without loss of generality,  $\mathcal A$ is represented by
\begin{equation}\label{matrix-representation}
\mathcal A=\left[
\begin{array}{cccc}
{0}_{n-2d}& & &  0\\
& \alpha_1 J & & \\
0 &   & \ddots \\
& & &  \alpha_d J
\end{array}
\right],\ \ \ J=\left[\begin{matrix}
0 & 1\\
-1& 0
\end{matrix}\right],
\end{equation}
see e.g. Rizzi \cite{Rizzi}. 
Here $0<\alpha_1\leq ...\leq \alpha_d,$ and $0_{n-2d}$ is the $(n-2d)\times (n-2d)$ square null-matrix. 
The layers $\mathfrak g_1$  and $\mathfrak g_2$ are generated by the left-invariant vector fields 
\begin{equation}\label{vector-field}
X_i=\partial_{x_i}-\frac{1}{2}\sum_{j=1}^k \mathcal A_{ij}x_j\partial_{x_{n+1}},\ \ i=1,...,n.
\end{equation}
Note that  $[X_i,X_j]=\mathcal A_{ij}\partial_{x_{n+1}},$  $i,j=1,...,n.$

If $n=2d$ (thus the kernel of $\mathcal A$  is trivial) and $\alpha_1=...=\alpha_d=4$, the Carnot group $\mathbb G$ reduces to the usual Heisenberg group $\mathbb H^d=\mathbb R^{2d}\times \mathbb R$. 

For our example, we shall consider a $6$-dimensional corank 1 Carnot group with the left-invariant vector fields given by  (\ref{vector-field}), by choosing $d=2$, $n=5$,  $\alpha_1=4$ and $\alpha_2=2.$ To be more explicit, the distribution $\mathcal D$ on  $(\mathbb G,\circ)$ is formed by the vector fields given by 
\begin{eqnarray} \label{DefIntMapH-1}
\left\{
\begin{array}{lll}
X_1=\partial_{x_1}; \\
X_2=\partial_{x_2}-2x_3\partial_{x_6};\\
X_3=\partial_{x_3}+2x_2\partial_{x_6};\\
X_4=\partial_{x_4}-x_5\partial_{x_6};\\
X_5=\partial_{x_5}+x_4\partial_{x_6}.\\
\end{array}\right.
\end{eqnarray}
Let $\textit{\textbf{a}}=(a_1,a_2,a_3,a_4,a_5)\in \Gamma(\mathcal D)$ given by the functions  
\begin{eqnarray} \label{DefIntMapH-2}
\left\{
\begin{array}{lll}
a_1=x_3^2x_5; \\
a_2=2{x_2}x_4x_6(x_6-2x_2x_3);\\
a_3=3x_1x_3x_5+4x_2^3x_4x_6;\\
a_4=x_2^2x_6(x_6-2x_4x_5);\\
a_5=x_1x_3^2+2x_2^2x_4^2x_6.\\
\end{array}\right.
\end{eqnarray}

\begin{example}\rm  We solve the problem 
	\begin{equation}\label{final-equation}
	X_iu=a_i\ {\rm in}\ \mathbb G,\ \ i=1,...,5.
	\end{equation}
\end{example}

To do this, we are going to fully explore Theorem \ref{fotetel}; by using the same notations, we identify  $A_1=0,$ $A_2=-2x_3$, $A_3=2x_2$, $A_4=-x_5$, $A_5=x_4.$ Moreover, $c_{23}=4=-c_{32}$, $c_{45}=2=-c_{54}$, the rest of the elements of the matrix $C=(c_{ij})$ being zero, $i,j=1,...,5$. In order to solve (\ref{final-equation}), we have to check relations (\ref{full-jellemzes}) and (\ref{full-jellemzes-2}), respectively.  It is easy to observe that (\ref{full-jellemzes}) is relevant only for $(i,j)=(2,3)$ and $(k,l)=(4,5)$ (the other choices giving always zero), where simple computations give that $X_2a_3-X_3a_2=8x_2^2x_4x_6$ and $X_4a_5-X_5a_4=4x_2^2x_4x_6$; thus, (\ref{full-jellemzes}) holds. Another simple reasoning shows that relation  (\ref{full-jellemzes-2}) is also verified; for instance,  $X_3X_2a_3-X_3X_3a_2=16x_2^3x_4=[X_2,X_3]a_3,$ the other relations following in the same way. 

Thus, Theorem \ref{fotetel} implies that the system (\ref{final-equation}) is solvable in $\mathcal F(\mathbb G)$; let $x_0={\bf 0}\in \mathbb G$ and any horizontal curve $\gamma=(\gamma_1,\gamma_2,\gamma_3,\gamma_4,\gamma_5,\gamma_6):[0,1]\to \mathbb G$ with $\gamma(0)={\bf 0}$ and $\gamma(1)=x=(x_1,x_2,x_3,x_4,x_5,x_6)\in \mathbb G$. Note that the horizontality of $\gamma$ means that
$$\dot\gamma_6=-2\gamma_3\dot\gamma_2+2\gamma_2\dot\gamma_3-\gamma_5\dot \gamma_4+\gamma_4\dot\gamma_5.$$
 Due to the latter relation and (\ref{C-V-1}),  some suitable rearrangements and $\gamma(0)={\bf 0}$ give that 
 \begin{eqnarray*}
 u(x)-c_0&=&\int_0^1 \sum_{i=1}^5a_i(\gamma(t))\dot \gamma_i(t)dt\\&=&\int_0^1 \frac{d}{dt}{ (\gamma_1(t)\gamma_3^2(t)\gamma_5(t))}dt+\int_0^1 \frac{d}{dt}{ (\gamma_2(t)\gamma_4(t)\gamma_6^2(t))}dt\\
 &=&\gamma_1(1)\gamma_3^2(1)\gamma_5(1)+\gamma_2(1)\gamma_4(1)\gamma_6^2(1)\\&=&
 x_1x_3^2x_5+x_2^2x_4x_6^2,
 \end{eqnarray*}
 for some $c_0\in \mathbb R$, which provides the solution of system  (\ref{final-equation}).

\section{Final remarks}\label{remark-section}

We conclude the paper with two remarks which can be considered as starting points of further investigations. \\

I)  \textit{Poincar\'e lemma on Finsler manifolds.} Let $(M,F)$ be an $m$-dimensional, not necessarily reversible Finsler manifold and $\Omega\subseteq M$ be a simply connected domain. Given a  vector field $\textit{\textbf{V}}\in \textit{\textbf{C}}^1(\Omega,T\Omega)$,   we are asking about the solvability of the equation  \begin{equation}\label{gradiens-megoldas-Finsler}
	{\bf \nabla}_F u=\textit{\textbf{V}}\ {\rm in}\ \Omega,
	\end{equation}
	where ${\bf \nabla}_F$ denotes the Finslerian gradient. Here, as usual ${\bf \nabla}_Fu(x)=J^*(x,Du(x))$, where  $J^*:T^*M\to TM$ is the Legendre transform  associating to each
	element $\alpha\in T_x^*M$ the unique maximizer on $T_xM$ of the map	$y\mapsto \alpha(y)-\frac{1}{2}F^2(x,y)$ and $Du(x)\in T_x^*M$ is 
 the  derivative of $u$ at $x\in M,$ see Ohta and Sturm \cite{Ohta-Sturm}. Note that in
	general,
	$u\mapsto{\bf \nabla}_F u $ is not linear. In order to solve (\ref{gradiens-megoldas-Finsler}), a necessarily curl-vanishing condition can be formulated by using the inverse  Legendre transform $J=(J^*)^{-1}$ and fundamental form of the Finsler metric $F$. However, we cannot adapt the proof of Theorem \ref{prop-1} into the Finsler setting. Indeed, we recall that in the proof of Theorem \ref{prop-1} we explored the metric compatibility and torsion-freeness of the Levi-Civita
	connection with respect to the given Riemannian metric; as we know, such properties are not simultaneously valid on a generic Finsler manifold unless it is Riemannian. \\
	
II)  \textit{Saint-Venant lemma on sub-Riemannian structures.} For simplicity, we shall consider only the usual  Heisenberg group $(\mathbb H^d,\mathcal D,g)$, where $\mathcal D=\{X_1,...,X_{2d}\}$ with $$X_{2i-1}=\partial_{x_{2i-1}}-2x_{2i}\partial_{x_{2d+1}}\ {\rm and}\  X_{2i}=\partial_{x_{2i}}+2x_{2i-1}\partial_{x_{2d+1}}, i=1,...,d,$$  and $g$ is the natural Riemannian metric on $\mathcal D,$ see (\ref{vector-field}.)  Given  a symmetric vector field $\textit{\textbf{e}}=(\textit{\textbf{e}}_1,...,\textit{\textbf{e}}_{2d})\in \Gamma(\mathcal D)^{2d}$ on $\Omega\subseteq \mathbb H^d$, i.e., ${e}_{ij}=e_{ji}$ for every $i,j=1,...,2d$ where $\textit{\textbf{e}}_i=\sum_{j=1}^{2d}e_{ij}X_j$, the question concerns the solvability of the  sub-Riemannian system 
	\begin{equation}\label{gradiens-megoldas-Finsler-szim}
	{\bf \nabla}_{s,H} \textit{\textbf{V}}=\textit{\textbf{e}}\ {\rm in}\ \Omega,
	\end{equation}
for the unknown vector field $\textit{\textbf{V}}=(V_1,...,V_{2d})\in \textit{\textbf{C}}^\infty(\Omega,\mathbb R^{2d}),$ 	where the components of the symmertric horizontal gradient ${\bf \nabla}_{s,H} $ are given by
	$$\frac{1}{2}(X_iV_k+X_kV_i),\ \ i,k=1,...,2d.$$
 The first challenging problem is to establish the necessary Saint-Venant compatibility relations associated to problem (\ref{gradiens-megoldas-Finsler-szim}) and then to apply Proposition \ref{Saint-Venant}; note that Schwartz type properties are not valid in this setting since usually $X_iX_j\neq X_jX_i$ for $i\neq j.$ Moreover, weaker versions of the Saint-Venant lemma on $\mathbb H^d$ would provide a sub-Riemannian Korn-type inequality as well.  Clearly, more general sub-Riemannian structures can also be considered instead of  Heisenberg groups verifying the assumptions of Theorem \ref{fotetel}. \\
	

\section{Appendix: proof of the Saint-Venant lemma (Proposition \ref{Saint-Venant})}\label{section-appendix}

%

A direct computation shows that if (\ref{szim-deriv}) has a solution, then the Saint-Venant 
compatibility relations (\ref{SV-relation}) trivially hold. 

Conversely, the  
Saint-Venant compatibility relations (\ref{SV-relation}) can be
written into the form
$${\partial}_{x_l}\left({\partial_{ x_j}e_{ik}}{}-{\partial_{ x_i}e_{jk}}\right)=
{\partial}_{x_k}\left({{\partial_
		{x_j}}e_{il}}-{{\partial_
		{x_i}}e_{jl}}\right),$$ which is
equivalent to
\begin{equation}\label{SV-cond}
{\partial}_{x_l}\left(g_{kt}g^{ts}\left({\partial_{ x_j}e_{is}}{}-{\partial_{ x_i}e_{js}}\right)\right)=
{\partial}_{x_k}\left(g_{lt}g^{ts}\left({{\partial_
		{x_j}}e_{is}}-{{\partial_
		{x_i}}e_{js}}\right)\right).
\end{equation}
If
$\textit{\textbf{W}}_{ij}$ is a vector field on $\Omega$ with the representation
$$\textit{\textbf{W}}_{ij}=W_{ijt}{\partial_{x_t}}=g^{ts}\left({\partial_{ x_j}e_{is}}{}-{\partial_{ x_i}e_{js}}\right){\partial}_{x_t},$$  relation (\ref{SV-cond}) can
be written equivalently into the form 
$${\partial}_{x_l}\left(g_{kt}W_{ijt}\right)=
{\partial}_{x_k}\left(g_{lt}W_{ijt}\right).$$ Thus, we
may apply Theorem 
\ref{prop-1}, i.e., there exists $p_{ij}\in C^2({\Omega})$ such that
$$\nabla_g p_{ij}=\textit{\textbf{W}}_{ij}\ \ {\rm on}\ \ \Omega,\ \ \forall i,j=1,...,m.$$
By components, the latter relation means that
$$g^{ts}{\partial_{x_s} p_{ij}}=W_{ijt}=g^{ts}\left({\partial_{ x_j}e_{is}}{}-{\partial_{ x_i}e_{js}}\right).$$
Multiplying from left by $g_{tl}$ and adding them, we have that
\begin{equation}\label{fontos-ez-p-ben}
{\partial_{x_l} p_{ij}}={\partial_{
		x_j}e_{il}}-{\partial_{ x_i} e_{jl}},\ \ \forall i,j,l=1,...,n.
\end{equation}
Since ${\partial_ {x_l}p_{ij}}+{\partial_{x_l}
	p_{ji}}=0,$ we can assume without loss of generality that
$p_{ij}+p_{ji}=0.$

If $q_{ij}=p_{ij}+e_{ij}$, then by $(\ref{fontos-ez-p-ben})$ we have
that
$${\partial_{x_k} q_{ij}}={\partial_{x_k} p_{ij}}+{\partial_{x_k} e_{ij}}={\partial_{x_j} e_{ik}}-{\partial_{x_i} e_{jk}}+{\partial_{x_k}
	e_{ij}}={\partial_{x_j} e_{ik}}+{\partial_{x_j} p_{ik}}={\partial_{x_j}
	q_{ik}}.$$ Again, the latter relation can be
transformed into
$${\partial _{x_k}}(g_{tj}g^{ts}q_{is})={\partial _{x_j}}(g_{tk}g^{ts}q_{is}).$$
Therefore, if $$\textit{\textbf{U}}_i=U_{il}{\partial}_{x_l}=g^{ls}q_{is}{\partial}_{x_l},$$ Theorem 
\ref{prop-1} implies the existence of $u_i\in C^2(\Omega)$ such
that $$\nabla_g u_i=\textit{\textbf{U}}_i,\ \ \forall i=1,...,m.$$ If we write the
components of the latter relation, it yields that
\begin{equation}\label{utolso-osszefugges}
{\partial_{x_l}u_i}=q_{il}, ,\ \ \forall i,l=1,...,m.
\end{equation}
Let $\textit{\textbf{V}}=(V_1,...,V_m)$  with 
%
$V_k=g^{ks}u_s,\ \ k=1,...,m.$ Consequently, by (\ref{utolso-osszefugges}), we have
\begin{eqnarray*}
	\frac{1}{2}\left({\partial_{x_i}}(g_{jk}V_k)+{\partial_{x_j}}(g_{ik}V_k)\right) &=&
	\frac{1}{2}\left({\partial_{x_i}}(g_{jk}g^{ks}u_s)+{\partial_{x_j}}(g_{ik}g^{ks}u_s)\right) \\
	&=& \frac{1}{2}\left({\partial_{x_i}}u_j+{\partial_{x_i}}u_i\right)  = \frac{1}{2}(q_{ij}+q_{ji}) \\&=&e_{ij},
\end{eqnarray*}
which is nothing but $\nabla_{s,g} \textit{\textbf{V}}=\textit{\textbf{e}}$, i.e.,  relation (\ref{szim-deriv}).
The Ces\`aro-Volterra integral formula follows at once by combining the above steps. \hfill $\square$
\\

 \noindent {\bf Acknowledgments.} The author thanks Professor
Philippe G. Ciarlet for his invitation  to the City University of
Hong Kong where the present work has been initiated. He is also grateful to Professors Ovidiu Calin and Der-Chen Chang for their suggestions and remarks.

\end{document}